\documentclass[a4paper,11pt]{book}

\usepackage[utf8]{inputenc}
\usepackage{geometry}
\usepackage[hidelinks]{hyperref}
\usepackage{titlesec}
\usepackage{etoolbox}
\usepackage{chngcntr}
\usepackage{graphicx}
\usepackage{ragged2e}
\usepackage{xcolor}
\usepackage{amsmath,amssymb,amsfonts,indentfirst}
\usepackage{fancyhdr}

\counterwithout{figure}{chapter}
\counterwithout{table}{chapter}
\counterwithout{equation}{chapter}

\patchcmd{\thebibliography}{\chapter*}{\section*}{}{}

\titleformat{\chapter}{\normalfont}{}{10pt}{\LARGE\bfseries\filcenter}
\newtheorem{teorema}{Theorem}[section]

\newtheorem{lema}{Lemma}[section]

\begin{document}
	
	\chapter{Pythagorean conics and homothetic triangles}
	
	\thispagestyle{empty} \vspace{-0.5cm} 
	
\begin{center}	
	Antonio A. Arcos Álvarez,  Emilio González Abril,  María-Jesús Vázquez-Gallo
\end{center}

{\bf Abstract}

This study investigates a generalisation of the Pythagorean theorem to the lengths of conic arcs constructed symmetrically on the sides of a right triangle. It is demonstrated that the theorem remains valid whenever the conic’s eccentricity is fixed and the ratio between the length of each arc’s sagitta and its corresponding side is constant. We identify the existence of a “Pythagorean centre”, defined as the common centre of all homotheties between the original right triangle and the enveloping right triangles of the infinite set of ``Pythagorean triples" of conics that can be derived from it. The proofs rely on the application of both geometrical and analytical techniques starting from classical Pythagoras theorem.

	\section{Introduction}

	This work proposes an extension of the Pythagorean theorem to arcs of conic curves which, under certain conditions, are supported on the sides of a right triangle.
	
	For the ancient Greeks, the Pythagorean theorem concerned the relationship between the areas of squares constructed on the sides of a right triangle \cite{Maor}. Generalisations of the theorem have been offered for other polygons, showing that the area of a polygon constructed on the hypotenuse of a right triangle is equal to the sum of the areas of similar polygons constructed on the triangle’s legs \cite{Bernal}. This last result was also recorded in Euclid’s Elements, but only for convex polygons \cite{PutzSpika}. However, as far as we know, no generalisations of the Pythagorean theorem have been established referring to the square of the length of curves constructed on the sides of a right triangle \cite{Agarwal}.
	
	This study demonstrates how, for arcs of conics with fixed eccentricity that rest at the endpoints of each side of a right triangle and are symmetric with respect to the perpendicular bisector of that side, the lengths of the respective conic arcs satisfy a quadratic relationship analogous to the Pythagorean theorem for the side lengths of the original right triangle, provided that the sagitta of each arc mantains the same constant proportional relationship with the side on which they rest.
	
	Having verified numerically that the result holds for several examples of conics (see \cite{sitioweb}), the deductive-axiomatic proof is handled both analytically and constructively, and several geometric observations of interest are presented in the process.
	
	This process begins by verifying that for each side of the triangle, there exists a unique non-degenerate conic arc of fixed eccentricity that touches the endpoints of the side and is symmetric with respect to the side’s perpendicular bisector (see Lemma \ref{unicidad}).
	
	It is then observed that, given a constant side–sagitta ratio for all three sides of the right triangle, the angle subtended by each side from the centre or the focus of the conic is the same for all three sides (see Lemma \ref{angulo}).
	
	Throughout the path to the final theorem (see Theorem \ref{teorema}), similarity and homothety are continuously present. Specifically, a homothetic relationship exists between the enveloping triangles of the various “Pythagorean triples” of conics that can be constructed on a single right triangle (see Figure \ref{figura8}). This homothety is centred at a singular point or “Pythagorean centre” which allows new Pythagorean triples of conics to be created via synthetic geometry, without the need of analytic calculation.
	
	From a mathematical perspective, the proof of the final result applies the classical arc length formula in terms of arc differentials in polar coordinates, together with the equality of the angles associated with each conic arc within the triple under consideration. This approach circumvents the difficulty that the integrals involved in the arc length of general ellipses and hyperbolas cannot be calculated explicitly, since the integrand functions do not possess elementary antiderivatives.
	
	\section{Results}
	
\begin{figure}[h!]
		\centering
		\includegraphics[width=0.4\textwidth]{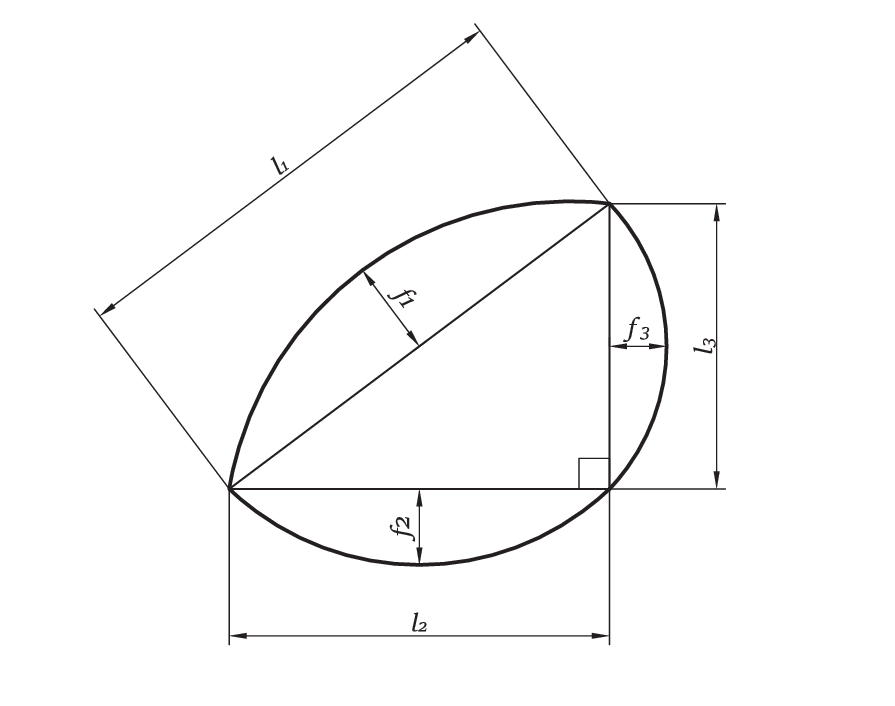}
		\caption{Conic arc constructed symmetrically on each side of a right triangle.}
		\label{figura1}
	\end{figure}

	Given a right triangle, on each of its sides we consider a non-degenerate conic arc of fixed eccentricity whose chord is the side in question, which is symmetric with respect to that side’s perpendicular bisector and for which the lengths of the corresponding sagittae are known (see Figure \ref{figura1}).
	
	This construction yields, for each side of the considered right triangle, three points through which the corresponding conic arc passes: two are the endpoints of the side, symmetric with respect to an axis of symmetry of the conic, and the remaining point, corresponding to the sagitta, belongs to that axis of symmetry. Let's verify that, under these conditions, for every side of the right triangle, the conic arc is unique.

	\begin{lema} $\label{unicidad}$
		
		Given a segment of length $l \in \mathbb{R}$, there exists a unique non-degenerate conic arc of fixed eccentricity $e \in [0,+\infty)$ constructed symmetrically on the segment.
		
	\end{lema}
	
	\textbf{Proof}
	
\begin{figure}[h!]
		\centering
		\includegraphics[width=0.35\textwidth]{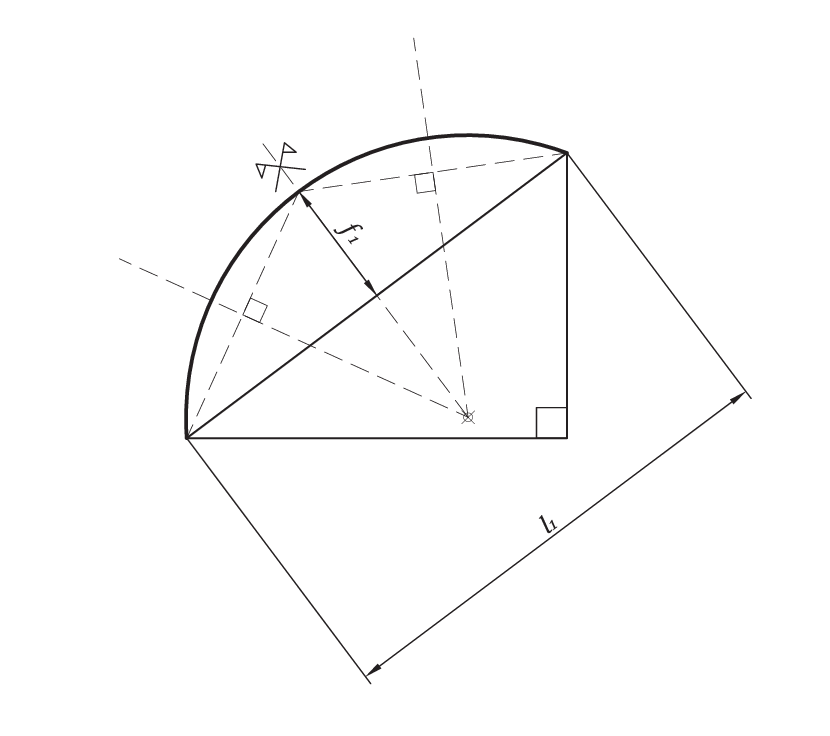}
		\caption{Uniqueness of the circular arc under the specified conditions.}
		\label{figura2}
	\end{figure}

	If the arc constructed on the segment is a circle (conic of eccentricity $e=0$), its uniqueness is clear: its centre must be the intersection point of the perpendicular bisectors of any pair chosen among the endpoints of the side of the right triangle and the point associated with the sagitta (see Figure \ref{figura2}).
	
\begin{figure}[h!]
		\centering
		\includegraphics[width=0.6\textwidth]{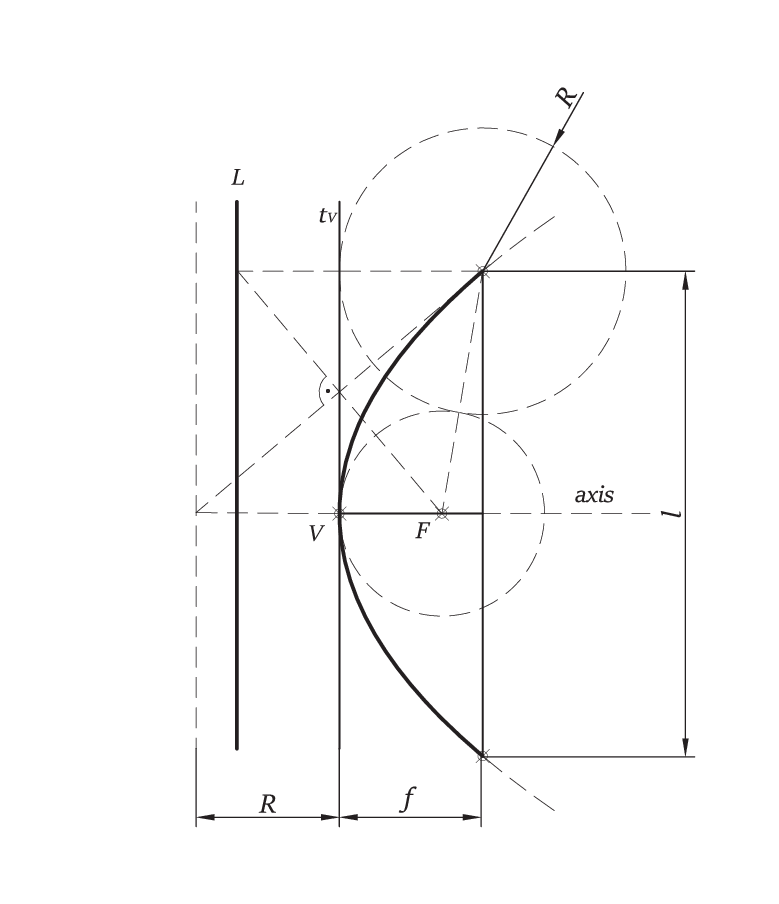}
		\caption{Construction of the unique parabola under the given conditions.}
		\label{figura3}
	\end{figure}

	For parabolas (conics of eccentricity $e=1$), uniqueness can also be demonstrated easily and constructively as follows (see Figure \ref{figura3}): The focus of the parabola may be found as the point on the axis of symmetry of the curve which is the centre of the circle tangent to the circle whose centre is one of the curve’s known points (the endpoints of the side of the triangle on which the curve is supported) and which is tangent to the tangent at the vertex of the parabola.
	
	The tangent at the vertex is simply the line perpendicular to the sagitta through its endpoint. Once the focus of the parabola is known, its directrix can be determined immediately and the parabola is therefore uniquely specified.
	
\begin{figure}[h!]
		\centering
		\includegraphics[width=0.5\textwidth]{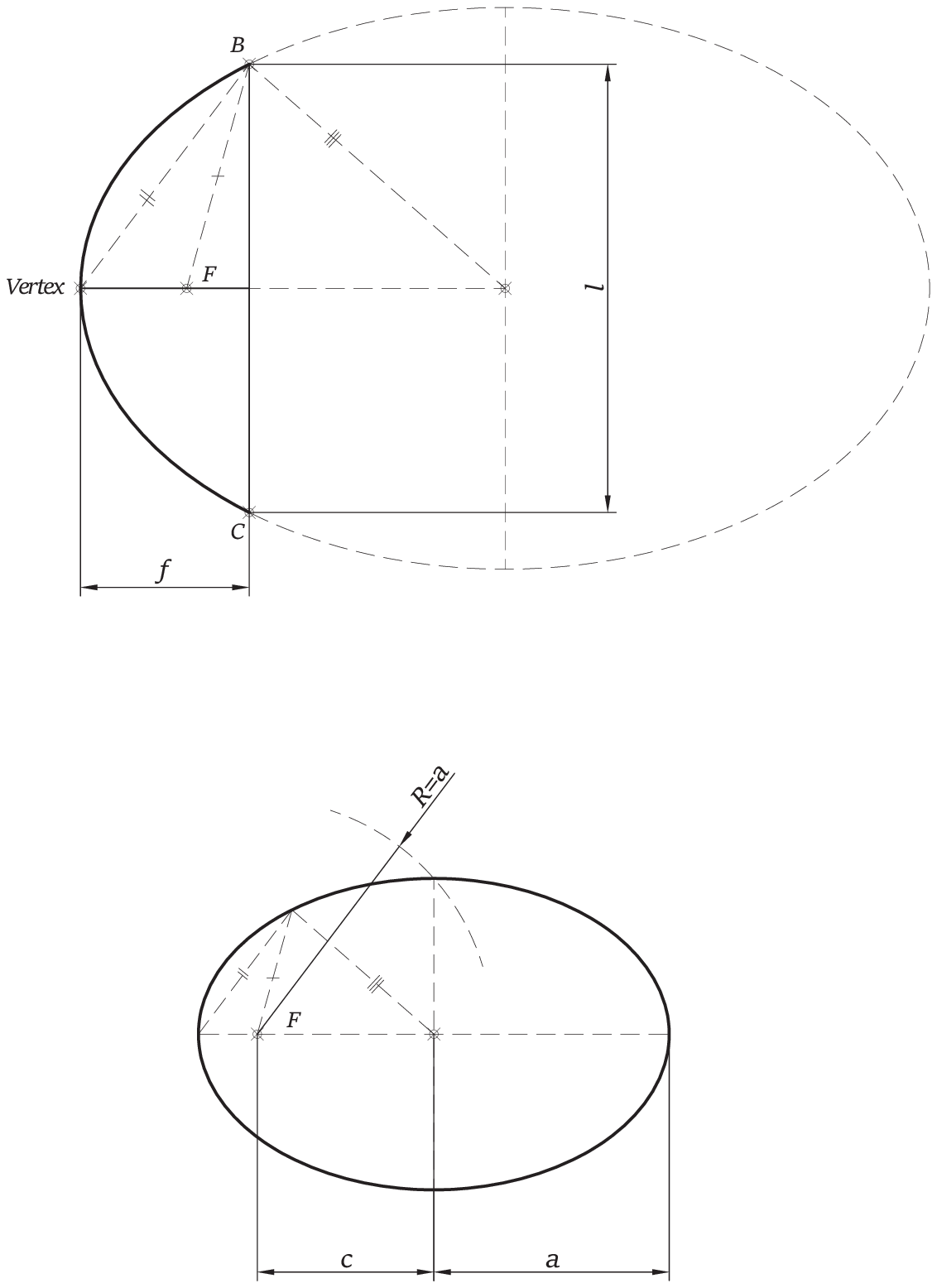}
		\caption{Elliptic arc under the specified conditions.}
		\label{figura4}
	\end{figure}

	In general, and constructively, for ellipses ($e<1$) and hyperbolas ($e>1$), the uniqueness of the considered conic arc is clear, noting that any two conics of equal eccentricity are similar (see Figure \ref{figura4} for the elliptical case). This concludes the proof of the lemma.

\
	
	The above result may also be proved analytically. For circles, it is known that a unique circle passes through any three non-collinear points. For parabolas, those points correspond to the graph of a quadratic polynomial determined uniquely by the data: the length of the triangle’s side $l$ and the sagitta $f$.
	
	In the case of ellipses with fixed eccentricity $0 < e < 1$, to verify the uniqueness of the ellipse arc symmetrically supported on a certain side of the right triangle, with given length $l$, we consider without loss of generality - rotating and/or translating- that the centre of the ellipse is at the origin of the Cartesian coordinates $x$ and $y$, and that its axes of symmetry are parallel to the coordinate axes. Then, the equation of such an ellipse is written as:
	$$
	\frac{x^2}{a^2} + \frac{y^2}{b^2} = 1,
	$$
	where $a$ and $b$ are the semi-axis lengths of symmetry, and the fact of having fixed the eccentricity $e = \frac{c}{a}$, with $a^2 = b^2 + c^2$ and $c$ being the distance from the centre of the ellipse to each of their foci, allows expressing $b$ in terms of $a$. Specifically:
	$$
	b^2 = a^2 (1 - e^2).
	$$

\begin{figure}[h!]
		\centering
		\includegraphics[width=0.45\textwidth]{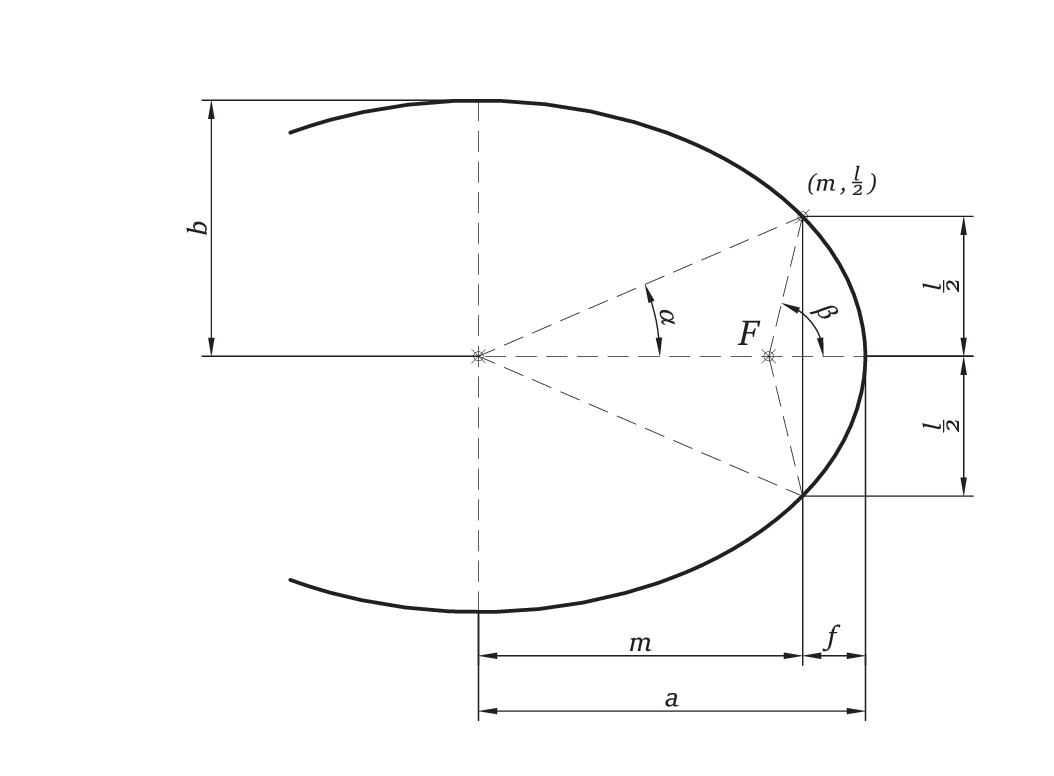}
		\caption{Elliptic arc symmetrically constructed on a segment of length $l$.}
		\label{figura5}
	\end{figure}

	In this situation, the equation of the ellipse must be satisfied for  three known points: the endpoints of the given length $l$ segment and one of the vertices of the ellipse. If this vertex is $(a,0)$, the other two points are $(m, \frac{l}{2})$ and its symmetric point $(m, -\frac{l}{2})$, for some $m > 0$ such that $a = m + f$, with $f$ being the corresponding sagitta length (see Figure \ref{figura5}). Then, a simple computation yields
	$$ m=\frac{l^2}{8f(1-e^2)}-\frac{f}{2}. $$
	The condition $m>0$ is equivalent to
	$\frac{b}{a}<\frac{\frac{l}{2}}{f}$, a requirement for constructing the ellipse.
	
	Thus, under the given conditions, the measures of the ellipse’s semi-axes
	$$ a=m+f= \frac{l^2}{8f(1-e^2)}+\frac{f}{2}, \quad b=\sqrt{a^2(1-e^2)}, $$
	are uniquely determined in terms of $e$, $l$, and $f$. The case where the known vertex is the point $(0,b)$ is similar. Therefore, the ellipse is unique.	
	
\begin{figure}[h!]
		\centering
		\includegraphics[width=0.45\textwidth]{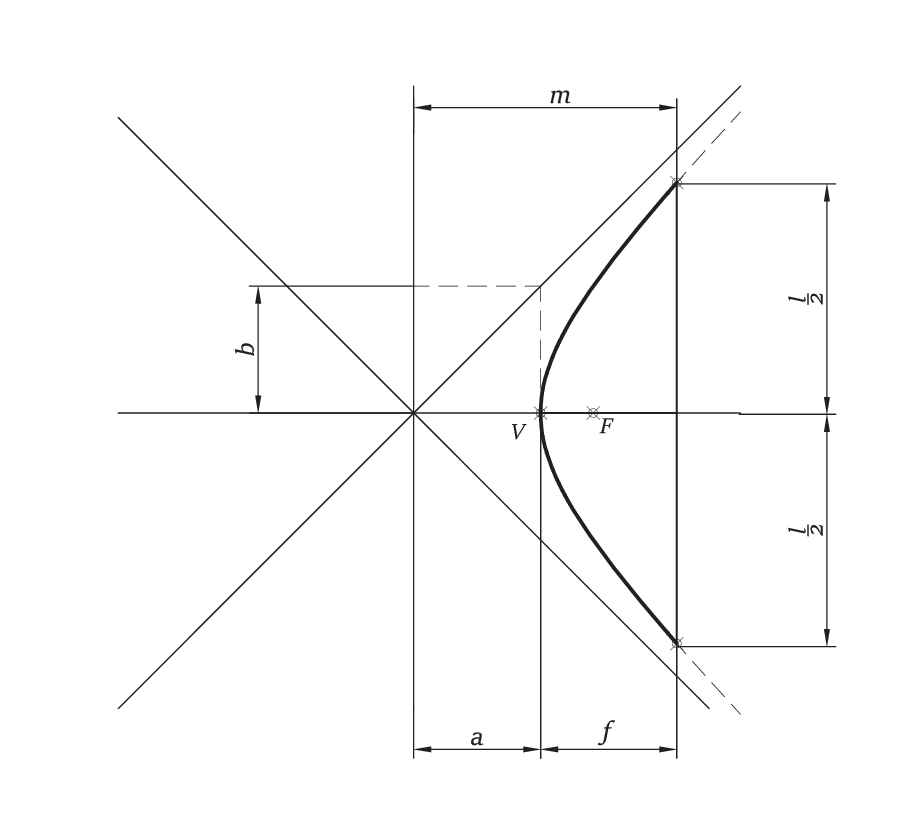}
		\caption{Hyperbolic arc constructed symmetrically on a segment of length $l$.}
		\label{figura6}
	\end{figure}

	Finally, the case of hyperbolas (conics with eccentricity $e>1$) is analogous to the case of ellipses, except that now the equation is
	$$\frac{x^2}{a^2}-\frac{y^2}{b^2}=1,$$
	with $a$ as the length of the semi-transverse axis, $a=m-f$, and $c^2=a^2+b^2$, where $c$ is the distance from focus to centre (see Figure \ref{figura6}). 
	
	Computing as before yields
	$$ m=\frac{l^2}{8f(e^2-1)}+\frac{f}{2}. $$
	Again, the lengths $a=m-f= \frac{l^2}{8f(e^2-1)}-\frac{f}{2}>0$ and $b=\sqrt{a^2(e^2-1)}$ are uniquely determined in terms of $e$, $l$, and $f<\frac{l}{\sqrt{e^2-1}}=\frac{l}{b}$ (to get $a>0$). Hence, the hyperbola is unique.
	
	\
	
	Once it is established that for each side of the initial right triangle there exists a unique conic arc under the given conditions, let us demonstrate that, with fixed side–sagitta ratio $k=\frac{l}{f}$, the angle subtended by each side of the triangle from the centre of the corresponding conic is the same for all three sides. The same holds for the angle from a focus of the conic. This will be useful in the extension of the Pythagorean theorem described in Theorem \ref{teorema}.
	
	\begin{lema} $\label{angulo}$
		
		Given a right triangle whose hypotenuse and legs have respective lengths $l_i \in \mathbb{R}$, with $i=1,2,3$, when considering the unique non-degenerate conic arc of fixed eccentricity $e \in [0,+\infty)$ constructed symmetrically on side $l_i$ with constant proportion $k:=\frac{l_i}{f_i}$, where $f_i \in \mathbb{R}$, $i=1,2,3$ is the length of each sagitta, it follows that the angle subtended by side $l_i$ from the centre of the respective conic, or from its focus, is constant.
		
	\end{lema}
	
	\textbf{Proof}
	
	The case of parabolas (eccentricity $1$) and circles (eccentricity $0$) are straightforward.
	
	For general ellipses (eccentricity $0<e<1$), as shown in the proof of Lemma \ref{unicidad}, if the known vertex is at the focal axis and we denote by $2\alpha$ the angle subtended by one of the triangle’s sides of length $l$ at the centre of the respective conic (see Figure \ref{figura4}), then the tangent of $\alpha$ is
	$$
	\tan \alpha = \frac{\frac{l}{2}}{m} = \frac{1 - e^2}{\frac{k^2 - 4(1 - e^2)}{4k}},
	$$
	which depends solely on the eccentricity $e$ and the fixed ratio $k$.
	
	Moreover, two angles with equal tangent differ by $\pi$ radians, but by construction $\alpha < \frac{\pi}{2}$, so the angle $\alpha$ is constant for all three sides of the triangle.
	
	Similarly, if we denote by $2\beta$ the angle subtended by one of the triangle’s sides of length $l$ at a focus of the respective conic (see Figure \ref{figura4}), then the tangent of $\beta$ is
	$$
	\tan \beta = \frac{\frac{l}{2}}{m - c} = \frac{\frac{l}{2}}{m - e a},
	$$
	which, using the expression for $m$ in the proof of Lemma \ref{unicidad}, also depends solely on the eccentricity $e$ and the fixed ratio $k$. If the known conic's vertex is not at the focal axis, the reasoning follows in a similar manner.
	
	The case of hyperbolae (eccentricity $e>1$) is also similar and it is left to the reader, thus concluding the proof.
	
	\
	
	From a constructive viewpoint, this result may again be proved using similarity between conics of fixed eccentricity.
	
	\
	
	Another interesting geometric observation, arising from the established conditions, is the exis\-tence of a homothety between the original right triangle and the enveloping triangles of the “Pythagorean conic” triples.
	
	\begin{figure}[h!]
		\centering
		\includegraphics[width=0.6\textwidth]{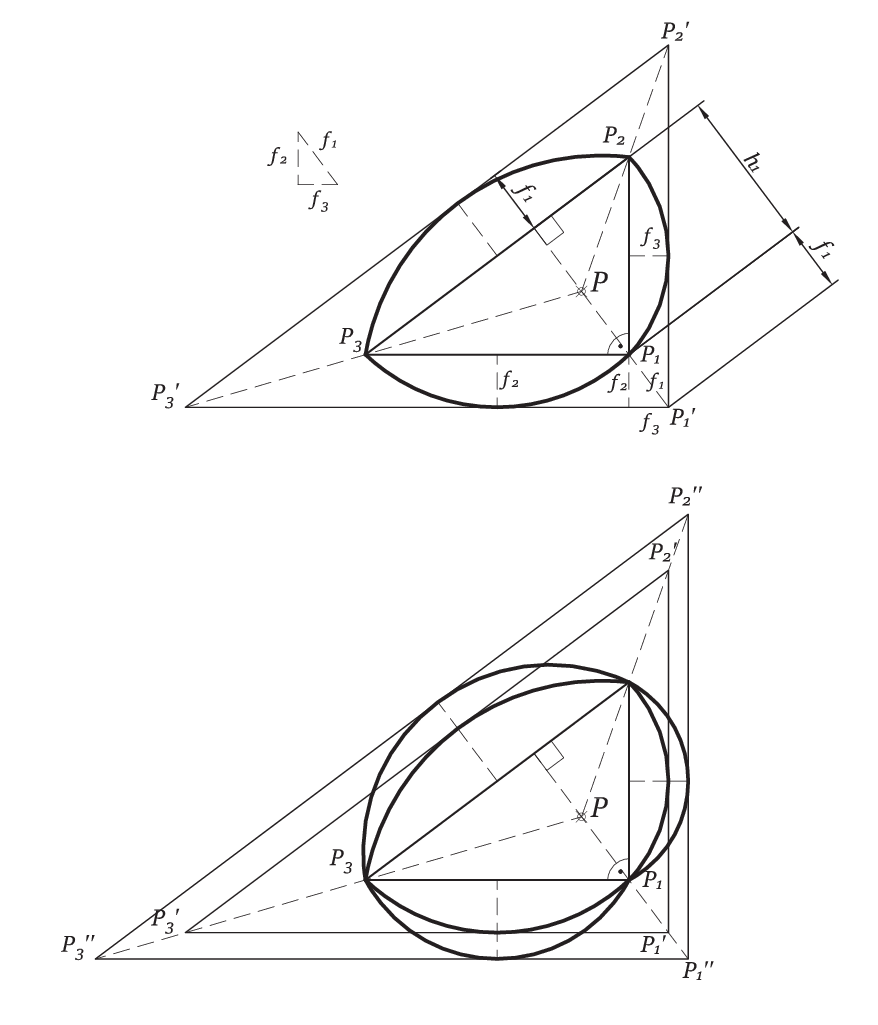}
		\caption{Homothetic triangles enveloping the triples of “Pythagorean conics”.}
		\label{figura7}
	\end{figure}
	
	Indeed, the imposed proportionality relation
	$$
	\frac{l_1}{f_1} = \frac{l_2}{f_2} = \frac{l_3}{f_3} = k,
	$$
	implies that the triangle whose sides are the sagittae $f_1$, $f_2$, and $f_3$ is similar to the original triangle, say $P_1 P_2 P_3$.
	
	The new enveloping triangle, say $P_1' P_2' P_3'$, constructed with sides parallel to the original through the endpoints of the sagittae $f_1$, $f_2$, and $f_3$, is also similar to the original.
	
	The segment $P_1 P_1'$ is therefore equal and parallel to $f_1$.
	
	If we consider the midpoint of the altitude $h_1$ corresponding to the right angle at $P_1$ and call it $P$, this point will always be the midpoint of the altitudes of the generated triangles $P_1' P_2' P_3'$, since its distance to $P_1'$ and to the opposite side $l_1'$ is always $\frac{h_1}{2} + f_1$.
	
	This condition implies that $P$ is the centre of homothety of the original triangle and all those generated as described before (see Figure \ref{figura7}).
	
	\
	
	\
	
	With all the above, an extension of the Pythagorean theorem may be demonstrated for conic arcs constructed on the sides of a right triangle under certain fixed conditions.
	
	\begin{teorema} \label{teorema}
		
		Given a right triangle with hypotenuse and legs of respective lengths $l_1, l_2, l_3$, let $c_i$, $i=1,2,3$, denote the length of the conic arc constructed symmetrically on side $l_i$ with fixed eccentricity $e \in [0,+\infty)$ and constant ratio $k := \frac{l_i}{f_i}$ for $i=1,2,3$, where $f_i \in \mathbb{R}$ is the length of each sagitta.
		
		Then:
		$$
		c_1^2 = c_2^2 + c_3^2.
		$$
		
	\end{teorema}
	
	\textbf{Proof:} Since the initial triangle is right-angled, if it holds that $c_i = g(k) l_i$ for $i=1,2,3$, with $k$ as the fixed ratio and $g(x)$ the same function in all three cases, then the result follows immediately, since
	$$
	c_2^2 + c_3^2 = g^2(k) l_2^2 + g^2(k) l_3^2 = g^2(k)(l_2^2 + l_3^2) = g^2(k) l_1^2 = c_1^2,
	$$
	where the classical Pythagorean theorem is applied in the penultimate equality.
	
	It suffices then to prove that, for $i=1,2,3$, the arc length $c_i$ can be expressed as the product of $l_i$, the length of the triangle side on which the arc rests, by a constant value $g(k)$, with $k$ the fixed ratio between the chord of the arc and its sagitta.
	
\begin{figure}[h!]
		\centering
		\includegraphics[width=0.3\textwidth]{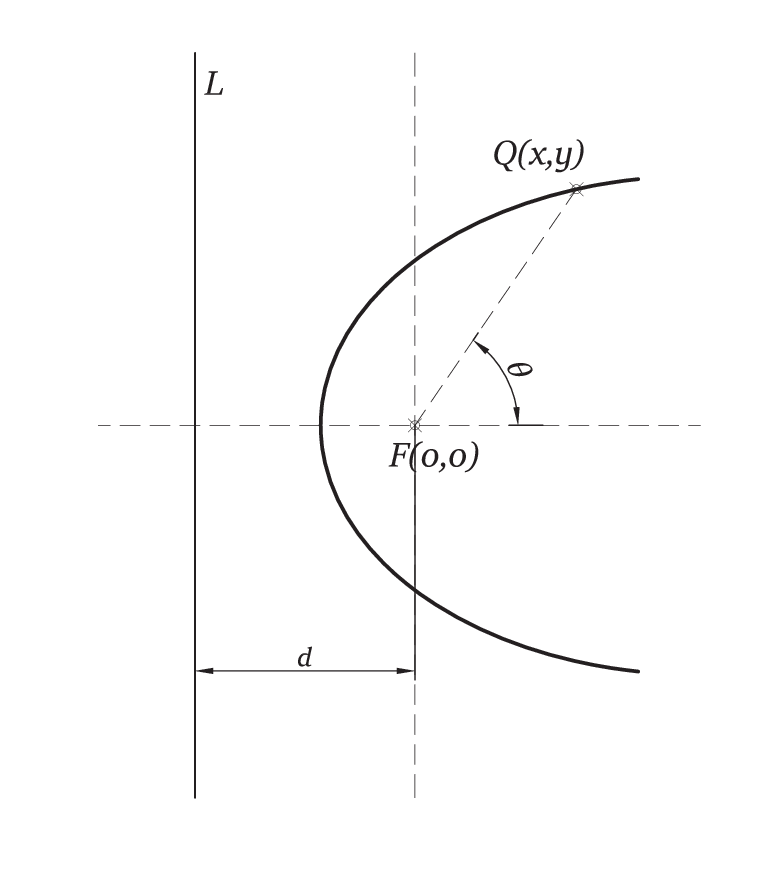}
		\caption{Conic arc defined by a line and a point exterior to it.}
		\label{figura8}
	\end{figure}

	On the one hand, it is known that in a plane, given a line $L$ and a point $F$ exterior to it, a conic of eccentricity $e>0$ is the locus of points $Q$ satisfying:
	$$
	d(Q,F)=e \cdot d(Q,L).
	$$

	where $d(Q,F)$ is the distance from $Q$ to the focus $F$ and $d(Q,L)$ is the distance from $Q$ to the directrix $L$. Developing the previous equality yields the polar equation of the conic:
	$$ r(\theta)=\frac{\pm d\,e}{1\mp e\cos \theta} $$
	where $d=d(F,L)$ denotes the distance between the focus $F$ and the directrix $L$, and for each point $Q$ on the conic, the angle $\theta$ is the one formed by the line joining $F$ to $Q$ and the polar axis of the conic (the axis of symmetry passing through $F$), while $r(\theta)$ is the distance from $Q$ to $F$, which we may—without loss of generality—set as the origin in Cartesian coordinates (see Figure \ref{figura8}).
	
	It is also known that the arc length of a planar curve given in polar coordinates by $r(\theta)$, between angle values $\theta_1$ and $\theta_2$, is given by the following integral, where $r'(\theta)$ is a continuous function.
	$$\int_{\theta_1}^{\theta_2} \sqrt{[r(\theta)]^2+[r'(\theta)]^2} \, d\theta.$$

	Due to the construction’s symmetry, in our case the interval of values for $\theta$ in this integral is also symmetric. Denoting this as $[-\alpha,\alpha]$, Lemma \ref{angulo} ensures $\alpha$ depends only on eccentricity $e$ and the fixed proportion $k=\frac{l_i}{f_i}$ for $i=1,2,3$. As for the integrand, calculation yields:
	$$\frac{de(1\pm 2e \cos \theta +e^2)}{(1\mp e \cos \theta)^4}.$$
	
	Hence,
	$$c_i=de \int_{-\alpha}^{\alpha} \frac{1\pm 2e \cos \theta +e^2}{(1\mp e \cos \theta)^4} d\theta,$$
	
	This is an elliptic integral that generally cannot be calculated explicitly, as the integrand does not possess an elementary antiderivative; however, here we use the fact that, regardless of its value, this integral depends only on the fixed eccentricity $e$ and the constant proportion $k$.
	
	Recalling that $d$ is the distance between the focus and the directrix of the corresponding conic, it follows that:
	\begin{itemize}
		\item If the conic is a parabola ($e=1$), then $d=de=\frac{kl_i}{8}=h(k)l_i$, with $h(x)=x/8$ for $i=1,2,3$.
		\item For ellipses with $0 < e < 1$, then $d=\frac{b^2}{c}$, with $b$ the length of the semi-axis perpendicular to the focal axis and $c$ the distance focus-centre. Consequently, $de=\frac{b^2}{a}=a(1-e^2)$, but from Lemma \ref{unicidad}, $a=\frac{kl_i}{8(1-e^2)}+\frac{l_i}{2k}=h(k)l_i$, with $h(x)=\frac{x}{8(1-e^2)}+\frac{1}{2x}$ for $i=1,2,3$.
		\item If the conic is a hyperbola ($e > 1$), again $d=\frac{b^2}{c}$ so $de=\frac{b^2}{a}=a(e^2-1)$, but as stated in Lemma \ref{unicidad}, $a=j(k)l_i$, with $j(x)=\frac{x}{8(1-e^2)}-\frac{1}{2x}$ for $i=1,2,3$.
	\end{itemize}
	This means that when the conic has eccentricity $e>0$, for $i=1,2,3$, it is indeed the case that $c_i=g(k)l_i$, with the same function $g(x)$, thus completing the proof in this case.
	
	For the remaining case, $e=0$, the conic is a special ellipse —a circle— its foci coincide with its centre, the lengths of its semi-axes equal its radius, and it is an exercise to verify that the theorem holds, which concludes the proof.
	
	\section{Conclusions}
	
	This work presents an extension of the Pythagorean theorem to non-degenerate conic arcs of fixed eccentricity, constructed symmetrically on the sides of a right triangle. Furthermore, a homothetic relationship is observed between the initial right triangle and the enveloping triangles of the infinite “Pythagorean triples” of conics that may be constructed on it. The centre of all these homotheties, or the “Pythagorean centre”, is always the same: the midpoint of the altitude of the relevant triangle with foot on the hypotenuse.
	
	This study suggests new questions regarding possible extensions of the Pythagorean theorem to arcs of other curves constructed on the sides of a right triangle, as well as possible variants of the presented results for geometries other than Euclidean, or possible practical applications of the geometric observations made.
	
	\begin{small}
		
	\end{small}

\end{document}